\documentclass[12pt, leqno,twoside]{article}
\usepackage{amssymb}
\usepackage{amsmath}
\usepackage{amsthm}
\usepackage{graphicx}
\usepackage[pagebackref,colorlinks,citecolor=blue,linkcolor=blue]{hyperref}
\usepackage{cite}   

\usepackage{ulem}
\usepackage{amsmath}
\usepackage{amsfonts}
\usepackage{amssymb}
\usepackage{amsmath,bbm,amssymb,amsxtra}
\usepackage{mathrsfs}
\usepackage{enumerate}
\usepackage{caption}
\allowdisplaybreaks
\usepackage{epsfig}
\usepackage{ae}
\def\vint{\mathop{\mathchoice%
         {\setbox0\hbox{$\displaystyle\intop$}\kern 0.22\wd0%
          \vcenter{\hrule width 0.6\wd0}\kern -0.82\wd0}%
         {\setbox0\hbox{$\textstyle\intop$}\kern 0.2\wd0%
          \vcenter{\hrule width 0.6\wd0}\kern -0.8\wd0}%
         {\setbox0\hbox{$\scriptstyle\intop$}\kern 0.2\wd0%
          \vcenter{\hrule width 0.6\wd0}\kern -0.8\wd0}%
         {\setbox0\hbox{$\scriptscriptstyle\intop$}\kern 0.2\wd0%
          \vcenter{\hrule width 0.6\wd0}\kern -0.8\wd0}}%
         \mathopen{}\int}

\newcommand{\bx}{{\partial \Omega}}

\newcommand{\loc}{\text{loc}}

\newcommand{\R}{{\mathbb R}}

\newcommand{\N}{{\mathbb N}}

\newcommand{\diam}{\text{\rm\,diam}}

\newcommand{\real}{{\mathbb R}}
\newcommand{\dist}{\text{\rm dist}}

\newcommand{\rarrow}{\rightarrow}

\newcommand{\mH}{{\mathcal H}}
\newtheorem{thm}{Theorem}[section]
\newtheorem{lem}[thm]{Lemma}
\newtheorem{prop}{Proposition}[section]
\newtheorem{rem}[thm]{Remark}
\newtheorem{defn}[thm]{Definition}%[section]
\newtheorem{example}[thm]{Example}
\numberwithin{equation}{section}

\makeatletter
\newcommand*{\mint}[1]{%
  \mint@l{#1}{}%
}
\newcommand*{\mint@l}[2]{%
  \@ifnextchar\limits{%
    \mint@l{#1}%
  }{%
    \@ifnextchar\nolimits{%
      \mint@l{#1}%
    }{%
      \@ifnextchar\displaylimits{%
        \mint@l{#1}%
      }{%
        \mint@s{#2}{#1}%
      }%
    }%
  }%
}
\newcommand*{\mint@s}[2]{%
  \@ifnextchar_{%
    \mint@sub{#1}{#2}%
  }{%
    \@ifnextchar^{%
      \mint@sup{#1}{#2}%
    }{%
      \mint@{#1}{#2}{}{}%
    }%
  }%
}
\def\mint@sub#1#2_#3{%
  \@ifnextchar^{%
    \mint@sub@sup{#1}{#2}{#3}%
  }{%
    \mint@{#1}{#2}{#3}{}%
  }%
}
\def\mint@sup#1#2^#3{%
  \@ifnextchar_{%
    \mint@sup@sub{#1}{#2}{#3}%
  }{%
    \mint@{#1}{#2}{}{#3}%
  }%
}
\def\mint@sub@sup#1#2#3^#4{%
  \mint@{#1}{#2}{#3}{#4}%
}
\def\mint@sup@sub#1#2#3_#4{%
  \mint@{#1}{#2}{#4}{#3}%
}
\newcommand*{\mint@}[4]{%
  \mathop{}%
  \mkern-\thinmuskip
  \mathchoice{%
    \mint@@{#1}{#2}{#3}{#4}%
        \displaystyle\textstyle\scriptstyle
  }{%
    \mint@@{#1}{#2}{#3}{#4}%
        \textstyle\scriptstyle\scriptstyle
  }{%
    \mint@@{#1}{#2}{#3}{#4}%
        \scriptstyle\scriptscriptstyle\scriptscriptstyle
  }{%
    \mint@@{#1}{#2}{#3}{#4}%
        \scriptscriptstyle\scriptscriptstyle\scriptscriptstyle
  }%
  \mkern-\thinmuskip
  \int#1%
  \ifx\\#3\\\else_{#3}\fi
  \ifx\\#4\\\else^{#4}\fi
}
\newcommand*{\mint@@}[7]{%
  \begingroup
    \sbox0{$#5\int\m@th$}%
    \sbox2{$#5\int_{}\m@th$}%
    \dimen2=\wd0 %
    \let\mint@limits=#1\relax
    \ifx\mint@limits\relax
      \sbox4{$#5\int_{\kern1sp}^{\kern1sp}\m@th$}%
      \ifdim\wd4>\wd2 %
        \let\mint@limits=\nolimits
      \else
        \let\mint@limits=\limits
      \fi
    \fi
    \ifx\mint@limits\displaylimits
      \ifx#5\displaystyle
        \let\mint@limits=\limits
      \fi
    \fi
    \ifx\mint@limits\limits
      \sbox0{$#7#3\m@th$}%
      \sbox2{$#7#4\m@th$}%
      \ifdim\wd0>\dimen2 %
        \dimen2=\wd0 %
      \fi
      \ifdim\wd2>\dimen2 %
        \dimen2=\wd2 %
      \fi
    \fi
    \rlap{%
      $#5%
        \vcenter{%
          \hbox to\dimen2{%
            \hss
            $#6{#2}\m@th$%
            \hss
          }%
        }%
      $%
    }%
  \endgroup
}

\begin{document} 
\title{\Large\bf $p$-harmonic mappings between metric spaces
\footnotetext{\hspace{-0.35cm}
 $2010$ {\it Mathematics Subject classfication}: 58E20, 46E35, 49Q10.
 \endgraf{{\it Key words and phases}: Metric valued Sobolev spaces, Dirichlet problem, upper gradients, Hajlasz-Sobolev space, trace operator.}
\endgraf{{\it ${}^{\mathbf{*}}$ Corresponding author}}
 }
}
\author{Chang-Yu Guo, Manzi Huang, Zhuang Wang$^*$ and Haiqing Xu}
\date{}
\maketitle

\begin{abstract}
	In this paper, we solve the Dirichlet problem for Sobolev maps between singular metric spaces that extends the corresponding result of Guo and Wenger [Comm. Anal. Geom. 2020]. The main new ingredient in our proofs is a suitable extension of the theory of trace for metric valued Sobolev maps developed by Korevaar and Schoen [Comm. Anal. Geom. 1993]. We also develop a theory of trace in the borderline case, which investigates a sharp condition to   characterize the existence of traces.
\end{abstract}

\tableofcontents

\section{Introduction}

%\subsection{Background and motivation}
The nonlinear Dirichlet problem associated to the $p$-harmonic mapping system in an Euclidean domain $\Omega\subset \R^n$, $1<p<\infty$, asks for a continuous map $u\colon \Omega\to \R^m$ so that 
\begin{equation*}
	\begin{cases}
		\nabla\cdot(|\nabla u|^{p-2}\nabla u)&=0 \quad \text{ in } \Omega,\\
	\qquad\qquad\quad\quad	u&=f                 \quad \text{ on } \bx.
	\end{cases}
\end{equation*}
The case $p=2$ corresponds to the classical Dirichlet boundary value problem associated to the harmonic mapping system. An equivalent way to formulate the general Dirichlet problem is to consider energy miniming mappings via the Euler-Lagrange equations. To be more precise, one considers minimizers of the $p$-energy 
$$E^p(u):=\int_{\Omega}|\nabla u|^p dx.$$

The above two formulations are not necessarily equivalent in general when we move from Euclidean spaces to Riemannian manifolds. Given two Riemannian manifolds $(M,g)$ and $(N,h)$, there is a natural $p$-energy functional acting on smooth maps
$$E^p(u):=\int_M |\nabla u|^p d\mu,$$
where $|\nabla u|$ is the Riemannian length of the gradient of $u$ and $\mu$ is Riemannian volume induced by $g$ on $M$.  Minimizers of the $p$-energy functional are called minimizing $p$-harmonic mappings, while ciritical points are called weakly $p$-harmonic mappings. In general, we only have the one-side inclusion: 
$$\{\text{minmizing }p \text{ harmonic mappings} \}\subset \{\text{weakly }p \text{ harmonic mappings} \}.$$ 
The case $p=2$ corresponds to the classical harmonic mappings. We refer the interested readers to \cite{sy97} for the theory of harmonic mappings, and to \cite{hl87,Luckhaus88,f90,gx19} for the theory of $p$-harmonic mappings, between Riemannian manifolds. 

One of the classical methods to solve the Dirichlet problem in the smooth setting is the direct method from the calculus of variations. To apply it, one essentially needs the following four ingredients:
\begin{itemize}
	\item A suitable $L^p$ theory for traces of manifold valued Sobolev maps;
	\item A global $p$-Poincar\'e inequality for Sobolev maps with zero trace;
	\item The Rellich-Kondrachov compactness theorem for manifold valued Sobolev maps;
	\item Lower semicontinuity of the energy functional $E^p$ with respect to $L^p$-convergence.
\end{itemize}
With all these ingredients at hand, the proof goes roughly as follows: Let $\{u_k\}\subset W^{1,p}(\Omega,N)$ be an energy minimizing sequence subordinate to the Dirichlet boundary condition $Tu_k=T\phi$ on $\bx$, where $\phi\in W^{1,p}(\Omega,N)$ is a fixed map and $T\colon W^{1,p}(\Omega,N)\to L^p(\bx,N)$ is the trace operator. The an easy application of the global Poincar\'e inequality, together with the characterization of traces of Sobolev maps, would give the boundedness of $\{u_k\}$ in $W^{1,p}(\Omega,N)$ (with respect to the Sobolev norm) and thus by the Rellich-Kondrachov compactness theorem for Sobolev spaces, we know that there exists a limiting map $u\in W^{1,p}(\Omega,N)$ such that a further subsequence $\{u_{k_i}\}$ of $\{u_{k}\}$ converges in $L^p(\Omega,N)$ to $u$. The lower semicontinuity of the energy functional implies the $p$-energy of $u$ would attain the minmimum, and at the same time, the convergence result for traces of Sobolev maps shall imply $Tu=T\phi$ on $\bx$. Therefore, $u$ is a proper solution to the Dirichlet problem.
 
Now, consider a mapping $u\colon X\to Y$, where $X=(X,d_X,\mu)$ is a metric measure space and $Y=(Y,d_Y)$ a metric space. Unlike the smooth Riemannian case, there is no natural $p$-energy functional associated to a sufficiently regular map. Indeed, there are several well-known (and generally different) $p$-energy functionals existing in the literature: the Korevaar-Schoen energy functional \cite{ks93}, the Jost energy functional \cite{j95}, the Hajlasz energy functional \cite{Ha96}, the upper gradient energy functional \cite{hk98,s00}, the Cheeger energy functional \cite{c99} and the Kuwae-Shioya energy functional \cite{ks03}; see \cite{Ha03,hkst12} for more energy functionals and the associated Sobolev spaces of metric valued maps. We would like to remark that the general interest in considering harmonic mappings in the singular metric setting dates back to the remarkable work of Gromov-Schoen \cite{gs92}, where the authors found important applications to rigidity problems for certain discrete groups; see \cite{g21} for a detailed survey on the theory of harmonic mappings between singular metric spaces.

In this article, we shall focus on the upper gradient energy functional and solve the associated Dirichlet problem. Throughout this paper, $X=(X,d_X,\mu)$ is assumed to be a complete metric measure space, $Y=(Y,d_Y)$ a complete   metric space, $\Omega\subset X$ a bounded domain and $\mH$ a $\sigma$-finite Borel regular measure on $\bx$. For notational simplicity, we sometimes drop the subscripts $X,Y$ from the distances $d_X, d_Y$ and simply write $d$.

%\subsection{Statement of the main results}
%Let $(X, d_X)$ and $(Y, d_Y)$ be two complete metric spaces, $\Omega\subset X$ a bounded domain and $\mu$ be a doubling measure on $\Omega$. We also fix a $\sigma$-finite Borel regular measure $\mH$ on $\bx$.

Before the statement of our main results, we recall a couple of definitions. One of the key concepts we shall need is the following definition of trace for metric valued functions.

\begin{defn}\label{trace-metric}
	Let $u\colon \Omega\rightarrow Y$ be a $\mu$-measurable function. Fix a point $x\in \bx$. If for some point $Tu(x)\in Y$, it holds
	\begin{equation}\label{trace-defn-Y}
		\lim_{r\rightarrow 0^+}\vint_{B(x, r)\cap\Omega}d_Y(u, Tu(x))\, d\mu=0, 
	\end{equation}
   then we say that the trace $Tu(x)$ of $u$ at $x\in \bx$ exists. Also, we say that $u$ has a trace $Tu$ on  $\bx$ if $Tu(x)$ exists for $\mH$-almost every $x\in \bx$.
\end{defn}
 
As in the smooth setting, we need to separate a class of admissible domains so that the Dirichlet problem is solvable. The following class of domains form a natural extension of the class of bounded Lipschitz domains in a smooth Riemannian manifold. 
%{\color{blue} check whether uniform domain is a special class. ZW: Not very sure, since it needs to determine the boundary measure $\mH$.}
\begin{defn}\label{def:p admissible domain}
We say that $\Omega\subset X$ is weakly $(q,\theta)$-admissible, $1<q<\infty$ and $\theta>0$, if 
\begin{itemize}
	\item $\mu$ is a doubling measure on $\Omega$;
    \item  $\mH$ is upper codimension-$\theta$ regular on $\partial\Omega$;
    \item $\Omega$ supports a local $q$-Poincar\'e inequality with $\theta<q<\infty$.
\end{itemize}
We say that $\Omega\subset X$ is $(q,\theta)$-admissible if in addition $\Omega$ supports a global $p$-Poincar\'e inequality for all $p\geq q$, that is, for $u\in N^{1,p}(\Omega)$ with $Tu=0$ $\mH$-almost everywhere on $\bx$, it holds
\begin{equation}\label{eq:Poincare for Sobolev functions with zero boundary value}
	\|u\|_{L^{p}(\Omega)}\leq C(\Omega)\|g_u\|_{L^p(\Omega)}.
\end{equation}
It is clear that if $\Omega$ is (weakly) $(q,\theta)$-admissible, then it is (weakly) $(p,\theta)$-admissible for any $p>q$.
\end{defn}

For the next concept, we refer to Section \ref{subsec:ultra-completion} below for the notion of a non-principal ultrafilter $\omega$ on $\mathbb{N}$ and the definition of ultra-limit $\lim_{\omega}a_m$ of a bounded sequence $\{a_m\}$ of real numbers. Let $(Y,d)$ be a metric space and $\omega$ a non-principal ultrafilter on $\mathbb{N}$. Denote by $Y_\omega$ the set of equivalent classes $[(y_m)]$ with the sequence $\{y_m\}$ in $Y$ satifying $\sup_md(y_1,y_m)<\infty$, where sequences $\{y_m\}$ and $\{y_m'\}$ are indentified if $\lim_{\omega}d(y_m,y_m')=0$. The metric space obtained by equipping $Y_\omega$ with the distance $d_{\omega}([(y_m)],[(y_m')])=\lim_{\omega}d(y_m,y_m')$ is called the ultra-completion or ultra-product of $Y$ with respect to $\omega$. It is clear that $Y$ isometrically embeds into $Y_\omega$ via the map $\iota\colon Y\to Y_\omega$, which assigns to $x$ the equivalent class $[(x)]$ of the constant sequence $\{x\}$. The following definition was introduced in \cite{gw17}.
\begin{defn}\label{def:1-complemented}
	A metric space $Y$ is said to be 1-complemented in some ultra-completion of $Y$ if there exists a non-principal ultrafilter $\omega$ on $\mathbb{N}$ for which there is a 1-Lipschitz retraction from $Y_\omega$ to $Y$.
\end{defn}
The class of metric spaces that are 1-complemented in some ultra-completion includes all proper metric spaces, all dual Banach spaces, some non-dual Banach spaces such as $L^1$, all Hadamard spaces and injective metric spaces; see \cite[Proposition 2.1]{gw17}.

Let $N^{1,p}(\Omega,Y)$ be the Sobolev space based on upper gradients and $E^p(u)$ the $p$-upper gradient energy functional of $u$ (see Section \ref{sec:preliminaries} below for precise definition). Our first main result can be formulated as follows.  
\begin{thm}\label{thm:main existence}
Suppose $\Omega\subset X$ is a $(q,\theta)$-admissible domain and $Y$ is a metric space that is 1-complemented in some ultra-completion of $Y$. Then for each $\phi\in N^{1,p}(\Omega,Y)$, $q<p$, there exists a mapping $u\in N^{1,p}(\Omega,Y)$ with $Tu=T\phi$ such that 
\begin{equation*}
	E^p(u)=\inf\left\{E^p(v): v\in N^{1,p}(\Omega,Y)\ \text{ and }\ Tv=T\phi\right\}.
\end{equation*}	
\end{thm}

To the best of our knowledge, Theorem \ref{thm:main existence} seems to be the most general setting for the solvability of the Dirichlet problem. In particular, it can be viewed as a natural extension of \cite[Theorem 5.6]{s01}, \cite[Theorem 1.4]{gw17} and \cite[Theorem 1.1]{g21}. 
In the formulation of Theorem \ref{thm:main existence}, we need the fact that the trace operator $T$ is well-defined on $N^{1,p}(\Omega,Y)$. When $\Omega\subset X$ is a bounded Lipschitz domain in a smooth Riemannian manifold and $Y$ is a complete metric space, this fact was established by Korevaar-Schoen in \cite[Section 12]{ks93}. 

Our second main result extends it to the more general singular setting.
\begin{thm}\label{thm-trace-1}
Suppose $\Omega$ is a weakly $(p,\theta)$-admissible domain and $Y$ is a complete metric space embedded isometrically into some Banach space. Then the trace operator 
$$T\colon N^{1, p}(\Omega, Y, d\mu)\rightarrow L^p(\bx, Y, d\mH)$$ 
is bounded and linear.
\end{thm}
Note that in Definition \ref{def:p admissible domain}, $p$ is assumed to be strictly larger than $\theta$ and in general Theorem \ref{thm-trace-1} fails for the borderline case $p=\theta$. In Section \ref{sec:refined trace}, we shall deduce sharper result in the borderline case by adding a weight $\omega$ to the underlying measure $\mu$; see Theorems \ref{thm-Banach-2} and \ref{thm-Banach-John} below.
 
Another crucial fact that we shall need in the proof of Theorem \ref{thm:main existence} is the following convergence result for traces of Sobolev maps with uniformly bounded energy.  When $\Omega\subset X$ is a bounded Lipschitz domain in a smooth Riemannian manifold and $Y$ is a complete metric space, this fact was established by Korevaar-Schoen in \cite[Theorem 1.12.2]{ks93}.
\begin{thm}\label{thm:convergence of trace}
Suppose $\Omega\subset X$ is weakly $(p,\theta)$-admissible and $Y$ is complete. Let $\{u_i\}\subset N^{1,p}(\Omega,Y)$ be a sequence with uniformly bounded energy, that is, 
$$\sup_{i\in \mathbb{N}}\int_{\Omega}g_{u_i}^pd\mu<\infty.$$
If $u_i$ converges to some $u\in N^{1,p}(\Omega,Y)$ in $L^p(\Omega,Y)$, then $Tu_i\to Tu$ in $L^p(\bx,Y)$.
\end{thm}

We next briefly comment on the ideas used in the proofs of our main theorems. As pointed out before, the proof of Theorem \ref{thm:main existence} relies essentially on the direct method from the calculus of variations. In the setting of Theorem \ref{thm:main existence}, a version of the Rellich-Kondrachov compactness theorem for metric valued Sobolev maps was obtained in \cite{gw17} and lower semicontinuity of the upper gradient energy is well-known, and thus the essential missing ingredient is a suitable $L^p$ theory for traces of metric valued Sobolev maps. 

The definition of trace in \cite[Section 1.12]{ks93} relies on the Lipschitz differentiable structure of $\Omega$, which looks apparently different than what we have introduced here. When the underlying spaces are metric measure spaces with much less geometric properties, Definition \ref{trace-metric} becomes a more natural way to define the trace.  When the target space $Y$ is $\real$, the trace results using Definition \ref{trace-metric} are under developing; see  \cite{KNW-agms,llw19,ls18,m17,ms18}. An useful observation in the proof of Theorem \ref{thm-trace-1} is that 
%we may embed $Y$ isometrically into a Banach space. Then 
by Lemma \ref{trace-isometric}, the isometric embedding (of $Y$ into some Banach space) and the trace operator commute, thus we only need to focus on the case when the target space $Y$ is a Banach space. 

In Theorem \ref{thm-trace-1}, it requires that $p>\theta$. It is natural to consider the borderline case when $\theta=p$. Theorem \ref{thm-Banach-2} and Theorem \ref{thm-Banach-John} deal with this borderline case. Especially  from Theorem \ref{thm-Banach-John} and Example \ref{exam:sharpness}, we obtain a sharp condition to full characterize  the existence of traces if additionally $\Omega$ is a John domain with compact closure.

%Comments on the proofs of Theorem 1.5 and 1.6, and then also the refiend theory of trace.

In this paper, we mainly consider the existence result for the Dirichlet problem. A natural question would be the interior regularity of the solutions. In the case when the target metric space $Y=\R$, there are local Lipschitz regularity results for solutions of the Dirichlet problem associated to the Cheeger energy and the upper gradient energy functional; see \cite{krs03,j14,bb11}. For general metric valued target space, there is a recent remarkable work due to Zhang and Zhu \cite{zz18}, where the authors derived local Lipschitz regularity of solutions of the Dirichlet problem associated to the Korevaar-Schoen energy functional; see also \cite{zzz19}. Recently, Guo and Xiang \cite{gx21} established local H\"older continuity of solutions of the Dirichlet problem associated to a variant of the Korevaar-Schoen $p$-energy functional. However, the method there relies crucially on the structure of the energy functional and does not extend to the upper gradient energy functional. We thus formulate it as an open question below.
\medskip

\textbf{Open question:} \textit{Under the assumptions of Theorem \ref{thm:main existence}, can we further establish local H\"older regularity result of the solution $u$? If so, under some kind of curvature assumption for $\Omega$ as in \cite{zz18} or \cite{krs03}, can we establish local Lipschitz regularity of $u$ for the harmonic case $p=2$?}
\medskip

The paper is organized as follows. In Section \ref{sec:preliminaries}, we recall the necessary definitions concerning metric valued Sobolev maps via upper gradients and ultra-completion of metric spaces. In Section \ref{sec:trace necessary}, we give an extension of the trace theory of Korevaar-Schoen and prove our trace theorem. Section \ref{sec:solution to Dirichlet} is devoted to the proof of Theorem \ref{thm:main existence}. In the final Section \ref{sec:refined trace}, we present a refined theory of trace in the borderline case.

\section{Preliminaries}\label{sec:preliminaries}
Let $(X, d_X,\mu)$ be a complete metric measure space and $(Y, d_Y)$ a complete   metric space. Let $\Omega\subset X$ be a bounded domain. We say that the measure $\mu$ is a {\it doubling measure} on $\Omega$ if there exists a constant $C_d\geq 1$ such that 
$$0<\mu(B(x, 2r)\cap\Omega) \leq C_d \mu(B(x, r)\cap \Omega)<\infty$$
for all $x\in \bar\Omega$ and $r>0$, where  $B(x, r):=\{y\in X: d(y, x)<r\}$ denotes an open ball centered at $x$ with radius $r$.

Given a set $F\subset \bar\Omega$ endowed with a $\sigma$-finite Borel regular measure $\mH$, we  say that $\mH$ is {\it upper codimension-$\theta$ regular on $F$} for some $\theta>0$ if there exists a constant $C_F$ such that
\begin{equation}\label{upper}
\mH(B(x, r)\cap F) \leq C_F \frac{\mu(B(x, r)\cap \Omega)}{r^\theta} 
\end{equation}
for  all $x\in F$ and $r>0$.

\subsection{Metric-valued Sobolev spaces via upper gradients}

Let $X=(X,d,\mu)$ be a metric measure space and $Z=(Z,d_Z)$ be a complete   metric space. Let $\Omega\subset X$ be a domain. 

For $p\geq 1$, we denote by $L^p(\Omega,Z)$ the space of all $\mu$-measurable and essentially separably valued map $u\colon \Omega\to Z$ such that for some $z_0\in Z$, the function $x\mapsto d(u(x),z_0)\in L^p(\Omega)$. A sequence $\{u_k\}\in L^p(\Omega,Z)$ is said to converge to $u\in L^p(\Omega,Z)$ if 
$$\int_{\Omega}d_Z^p(u(x),u_k(x))d\mu(x)\to 0 \qquad \text{as }k\to\infty.$$
When $(Z,d_Z)=(V, |\cdot|)$ is a Banach space, we may edow $L^p(\Omega,V)$ with a natural norm 
$$\|f\|_{L^p(\Omega, V)}:=\left(\int_\Omega |f|^p\, d\mu\right)^{1/p}.$$
Similarly, we can define $L^p(\bx, V):=L^p(\bx, V, d\mH)$. If $V$ is $\real$, we set $L^p(\Omega, \real)=:L^p(\Omega)$ and $L^p(\bx, d\mH)=:L^p(\bx)$ for brevity. 

We next introduce metric valued Sobolev spaces based on upper gradients. This concept was first introduced in~\cite{hk98} and then functions with $p$-integrable upper gradients were studied in~\cite{km98}. Later, the theory of real-valued Sobolev spaces based on upper gradients was explored in-depth in~\cite{s00}. Here we only give a very brief introduction and refer the interested readers to the recent monograph \cite{hkst12} for more information.

\begin{defn}[Upper gradients]\label{def:upper gradients} A Borel function $g\colon \Omega\rightarrow [0,\infty]$ is called an upper gradient for a map $u\colon \Omega\to Z$ if for every rectifiable curve $\gamma\colon [a,b]\to \Omega$, we have the inequality
\begin{equation}\label{upper-gradient}
	d_Z(u(\gamma(b)),u(\gamma(a)))\leq \int_\gamma g\ ds\text{.}
\end{equation}
If inequality \eqref{upper-gradient} holds for $p$-almost every curve, then $g$ is called a $p$-weak upper gradient for $u$.  
\end{defn}

A $p$-weak upper gradient $g$ of $u$ is minimal if for every $p$-weak upper gradient $\tilde{g}$ of $u$, $\tilde{g}\geq g$ $\mu$-almost everywhere.  If $u$ has an upper gradient in $L^p_{\loc}(\Omega)$, then $u$ has a unique (up to sets of $\mu$-measure zero) minimal $p$-weak upper gradient.  We denote the minimal upper gradient by $g_u$. The Sobolev space $N^{1,p}(\Omega,Z)$ consists of all $u\in L^p(\Omega,Z)$ with an $L^p$-integrable minimal $p$-weak upper gradient $g_u\in L^p(\Omega)$. For each $u\in N^{1,p}(\Omega,Z)$, we shall use $E^p(u)$ to denote the upper gradient energy functional of $u$, that is, 
$$E^p(u)=\int_{\Omega} g_u^pd\mu.$$

An alternative way to introduce $N^{1,p}(\Omega,Z)$ is to use isometric embedding $Z\subset V$ and then define $N^{1,p}(\Omega,Z)$ as the Banach space-valued Sobolev spaces $N^{1,p}(\Omega,V)$. As this will be convenient for us later in establishing the theory of trace, we briefly record Banach space valued Sobolev spaces $N^{1,p}(\Omega,V)$ here.
%A Borel function $g: \Omega\rightarrow [0, \infty]$ is said to be an {\it upper gradient} of a map $u: \Omega\rightarrow V$ if 
%\begin{equation}\label{upper-gradient}
%|u(\gamma(a))-u(\gamma(b))|\leq \int_{\gamma} g\, ds
%\end{equation}
%for every rectifiable curve $\gamma: [a, b]\rightarrow X$. We say that $g$ is a {\it $p$-weak upper gradient} of $u$ if \eqref{upper-gradient} holds for $p$-a.e. rectifiable curves in $\Omega$.

The {\it Dirichlet space} $D^{1, p}(\Omega, V)$ consists of all measurable functions $u\colon \Omega\rightarrow V$ that have an upper gradient belonging to $L^p(\Omega)$. We can equip the Dirichlet  space $D^{1, p}(\Omega, V)$ with the seminorm
$$\|u\|_{D^{1, p}(\Omega, V)}:=\inf_g\|g\|_{L^p(\Omega)},$$
where the infimum is taken over all $p$-weak upper gradient $g$ of $u$.

Let
$$\tilde N^{1, p}(\Omega, V) =D^{1, p}(\Omega, V) \cap L^p(\Omega, V)$$
be equipped with the seminorm
$$\|u\|_{\tilde N^{1, p}(\Omega, V)} =\|u\|_{L^p(\Omega, V)}+\|u\|_{D^{1, p}(\Omega, V)}.$$
We obtain a normed space $N^{1, p}(X, V)$, which is called the {\it Sobolev space} of $V$-valued functions on $\Omega$, by passing to equivalence classes of functions in $\tilde N^{1, p}(\Omega, V)$, where $u_1\sim u_2$ if and only if $\|u_1-u_2\|_{\tilde N^{1, p}(\Omega, V)}=0$. Thus,  
$$N^{1, p}(\Omega, V):=\tilde N^{1, p}(\Omega, V)/\{u\in \tilde N^{1, p}(\Omega, V): \|u\|_{\tilde N^{1, p}(\Omega, V)=0}\}.$$
Since we may embed the metric space $Z$ isometrically into some Banach space   $L^\infty(Y)$, we can alternatively define $N^{1,p}(\Omega,Z)$ via $N^{1,p}(\Omega,Z):=N^{1,p}(\Omega,L^\infty(Z))$; see \cite[Section 7]{hkst12}.

We say $\Omega$ supports a local $q$-Poincar\'e inequality, $1\leq q<\infty$, if there exist constants $C>0$ and $\lambda\geq 1$ such that
\begin{equation}\label{eq:local Poincare inequality}
\vint_{B(x, r)\cap \Omega} |u-u_{B(x, r)\cap \Omega}|\, d\mu\leq C r \left(\vint_{B(x, \lambda r)\cap\Omega}g^q\, d\mu\right)^{1/q}
\end{equation}
holds for all $x\in \Omega$ and $r>0$, for every function $u: \Omega\rightarrow \real$ belonging to $L^1_{\rm loc}(\Omega)$, and for every upper gradient $g$ of $u$.

\begin{rem}\rm
(i) If $\mu$ is a doubling measure on $\Omega$, then the inequality \eqref{eq:local Poincare inequality}  holds not only for  $x\in \Omega$, but also for $x\in \bx$; see \cite[Remark 2.13]{m17} for a discussion.

(ii) It follows from \cite[Theorem 8.1.42]{hk98} that for doubling metric measure spaces, the validity of a Poincar\'e inequality is independent of the target Banach spaces. Hence if $\mu$ is a doubling measure on $\Omega$ and $V$ is a Banach space, then that $\Omega$ supports a local $q$-Poincar\'e inequality implies that  \eqref{eq:local Poincare inequality} holds for for all $x\in \bar\Omega$ and $r>0$, for every function $u: \Omega\rightarrow V$ belonging to $N^{1, p}(\Omega, V)$ with every upper gradient $g$ of $u$.

(iii) If we instead consider metric space valued functions, by embedding the target metric space into a Banach space,  that $\Omega$ supports a local $q$-Poincar\'e inequality implies that there exist constants $C>0$ and $\lambda\geq 1$ such that
\[\vint_{B(x, r)\cap \Omega}\vint_{B(x, r)\cap \Omega}d_Z(u(z),u(y))\, d\mu(z)\, d\mu(y)\leq C r \left(\vint_{B(x, \lambda r)\cap\Omega}g^q\, d\mu\right)^{1/q}\]
holds for all $x\in \Omega$ and $r>0$, for every map $u\colon \Omega\rightarrow Z$ belonging to $L^1_{\rm loc}(\Omega, Z)$, and for every Borel function $g$ being upper gradients of $u$.
\end{rem}

%\subsection{A short discussion on the global Poincar\'e inequality for Sobolev functions with zero trace}

In Definition \ref{def:p admissible domain}, we imposed a global $p$-Poincar\'e inequality on $\Omega$ for Sobolev functions with zero trace. When $\Omega$ is a bounded Lipschitz domain in a Riemannian manifold, this condition is easily verified. When $\Omega$ is a bounded domain in a general metric measure space $X$, it seems to be not yet clear what should be a reasonable geometric assumption imposed on $\Omega$. One essentially needs that the zero extension of a Sobolev function on $\Omega$ with zero trace shall be a global Sobolev function on $X$. For this, we need to analyze the size of the exceptional set of points on $\partial\Omega$ such that \eqref{trace-defn-Y} fails. This involves a careful study of the pointwise behavior of a Sobolev function in singular metric spaces and it goes beyend the scope of this paper. We refer the interested readers to \cite{kkst12} for results and discussions along this direction.

%\begin{KET}
%Every metric space $Y$ embeds isometrically into the Banach space $L^\infty(Y)$.
%\end{KET}
%{\color{red}Definitions of  Metric-valued Sobolev space? Do we use $N^{1, p}(X, L^\infty)$?}

%Hence instead of studying the metric-valued Sobolev spaces and their traces, it is equivalent to consider the Banach space-valued Sobolev spaces and their traces.

\subsection{Ultra-completions of metric spaces}\label{subsec:ultra-completion}
We briefly recall the relevant definitions concerning ultra-completions and ultra-limits of metric spaces. Details can be found for instance in \cite{bh99}.

A non-principal ultrafilter on $\N$ is a finitely additive probability measure $\omega$ on $\N$ such that every subset of $\N$ is measurable and such that $\omega(A)$ equals 0 or 1 for all $A\subset \N$ and $\omega(A)=0$ whenever $A$ is finite. Given a compact Hausdorff topological space $(Z,\tau)$ and a sequence $\{z_m\}\subset Z$ there exists a unique point $z_\infty\in Z$ such that $\omega(\{m\in \N:z_m\in U\})=1$ for every {  $U\ni \tau$} containing $z_\infty$. We denote the point $z_\infty$ by $\lim_\omega z_m$.

Let $Y=(Y,d)$ be a metric space and $\omega$ a non-principal ultrafilter on $\N$. A sequence $\{y_m\}\subset Y$ is bounded if $\sup_m d(y_1,y_m)<\infty$. Define an equivalence relation $\sim$ on bounded sequences in $Y$ by considering $\{y_m\}$ and $\{y_m'\}$ equivalent if $\lim_\omega d(y_m,y_m')=0$. Denote by $[(y_m)]$ the equivalence class of $\{y_m\}$. The ultra-completion $Y_\omega$ of $Y$ with respect to $\omega$ is the metric space given by the set
$$Y_\omega:=\{[(y_m)]:\{y_m\} \text{ bounded sequence in }Y\},$$
equipped with the metric 
$$d_\omega([(y_m)],[(y_m')]):=\lim_\omega d(y_m,y_m').$$ 
The ultra-completion $Y_\omega$ of $Y$ is a complete metric space, even if $Y$ itself is not complete.

\section{Extension of the trace theory of Korevaar-Schoen}\label{sec:trace necessary}

%\subsection{Extension of the trace theory of Korevaar-Schoen}
Let $(X, d_X,\mu)$ be a complete metric measure space, $(V, |\cdot|)$ be a Banach space  and $\Omega\subset X$ be a bounded domain. Assume $\bx$ is endowed with  an upper codimension-$\theta$ regular measure $\mathcal H$ with $\theta>0$. 

We give an alternative definition of the trace for Banach valued maps.
\begin{defn}\label{trace-Banach}
Let $u\colon \Omega\rightarrow V$ be a $\mu$-measurable function. Then $Tu(x)\in V$ is the trace of $u$ at $x\in \bx$ if the following equation holds:
\begin{equation}\label{trace-defn}
\lim_{r\rightarrow 0^+}\vint_{B(x, r)\cap\Omega}|u-Tu(x)|\, d\mu=0.
\end{equation}
We say that $u$ has a trace $Tu$ on  $\bx$ if $Tu(x)$ exists for $\mH$-almost every $x\in \bx$.
\end{defn}

We next show that Definition \ref{trace-Banach} is consistent with Definition \ref{trace-defn-Y}.
\begin{lem}\label{trace-isometric}
Let $(Y, d_Y)$ be a complete metric space and $h\colon Y\rightarrow L^\infty(Y)$ be an isomeric embedding. For any $\mu$-measurable function $u\colon \Omega\rightarrow Y$, if the function $h\circ u\colon \Omega\rightarrow L^\infty(Y)$ has a trace in the sense of Definition \ref{trace-Banach}, then the function $u$ has a trace in the sense of Definition \ref{trace-metric}.
\end{lem}
\begin{proof}
Let $T(h\circ u)$ be the trace of $h\circ u$ in the sense of Definition \ref{trace-Banach}. Then for $\mH$-a.e. $x\in \bx$, the equation \eqref{trace-defn} holds for $T(h\circ u)$. Note that \eqref{trace-defn} implies that for each $k\in \N$, we can find a ball $B(x, r_k)$ centered at $x\in \bx$ with radius $r_k>0$ such that there exists a point $x_k\in B(x,r_k)\cap \Omega$ with 
$$|T(h\circ u)(x)-h\circ u(x_k)|<2^{-k}.$$ 
Hence $\{h\circ u (x_k)\}_{k\in \N}$ is a Cauchy sequence in $L^\infty(Y)$ that converges to $T(h\circ u)(x)$. 

Since $h$ is an isometric embedding, $\{u (x_k)\}_{k\in \N}$ is a Cauchy sequence in $Y$ and hence has a limit in the complete metric space $Y$, for which we denote by $Tu(x)$. Moreover, it follows from the isometric property of $h$ that $h\circ Tu(x)=T(h\circ u)(x)$ and \eqref{trace-defn-Y} is satisfied with  $T(h\circ u)$ for $\mH$-a.e. $x\in \bx$. Thus, $Tu$ is the trace of the function $u$ in the sense of Definition \ref{trace-metric}.
\end{proof}

%\begin{rem}\label{rmk:on trace equivalence}
The proof of Lemma \ref{trace-isometric} actually tells that the isometric embedding $h$ and the trace operator $T$ commute. Thus to develop a theory of trace, we shall not distinguish the traces operators in Definitions \ref{trace-metric} and \ref{trace-Banach}. From now on, we shall focus on the case when $(Y,d_Y)=(V,|\cdot|)$ is a Banach space.
%\end{rem}
%{\color{red} Here the above lemma actually tells that the isometric embedding and the trace operator commutes, so I think it is possible to write the commutation property as a lemma and we do not distinguish the traces operators in Definition \ref{trace-metric} and \ref{trace-Banach}.  Which one do you think better?}

For any $f\in L^p_{\rm loc}(\Omega, V)$, we define the centered fractional maximal operator as
\begin{equation}\label{maximal-operator}
M_{\theta, p} f(z)=\sup_{0<r<2\diam (\bx)} \left(r^\theta\vint_{B(z, r)\cap \Omega} |f|^p\, d\mu\right)^{1/p},\ \ \text{ for each } \ z\in \bx.
\end{equation}
Then it is easy to see that this fractional maximal operator maps $L^p_{\rm loc}(\Omega, V)$ into the space of real-valued lower semicontinuous functions on $\bx$.

The following boundedness result follows from  \cite[Lemma 4.2]{m17}.
\begin{lem}\label{fractional-maximal}
Let $1\leq p<\infty$. Then the fractional maximal operator $M_{\theta, p}$ is bounded from $L^p(\Omega)$ to weak-$L^p(\bx)$.
\end{lem}

We are ready to prove the boundedness of the trace operator for Banach space valued Sobolev maps. When the Banach space is $\real$, the result was obtained in \cite{m17}. The essential idea of the proof is similar with the one used in  \cite{m17}.
\begin{thm}\label{thm-Banach-1}
Suppose $\Omega\subset X$  is weakly $(p,\theta)$-admissible for some $p>1$. 
Then the trace operator $T\colon N^{1, p}(\Omega, V)\rightarrow L^p(\bx, V)$ is bounded and linear.
\end{thm}
\begin{proof}
Let $u\in N^{1, p}(\Omega, V)$ and $R=2\diam(\Omega)$ be fixed. For any $z\in \bx$ and $k\in \N$, we define 
$$T_k u(z)=\vint_{B(z, 2^{-k} R)\cap\Omega} u\, d\mu.$$ 

We first show that the limits
$$\widetilde Tu=\lim_{k\rightarrow \infty} T_k u$$
exist $\mH$-almost everywhere on $\bx$.
It suffices to show that the function 
$$u^*=\sum_{k\geq 0} |T_{k+1}u -T_k u|+|T_0 u|$$
belongs to $L^p(\bx)$, since $u^*\in L^p(\bx)$ implies that $u^*(z)<\infty$ for $\mH$-almost everywhere $z\in \bx$.
Then it suffices to show that
$$\|u^*\|_{L^p(\bx)}\leq \|T_0 u\|_{L^p(\bx, V)}+\sum_{k\geq 0}\|T_{k+1}u- T_k u\|_{L^p(\bx, V)}<\infty.$$

Notice that $T_0 u(z)=\vint_{B(z, R)\cap\Omega} u\, d\mu=\vint_{\Omega} u\, d\mu$ for any $z\in \bx$, since $\Omega\subset B(z, R)$ for any $z\in \bx$. It follows from the upper codimension relation \eqref{upper} that
\[\|T_0 u\|^p_{L^{p}(\bx, V)}\leq \int_{\bx} \vint_{B(z, R)\cap \Omega} |u|^p\, d\mu\, d\mH(z)\lesssim \int_{\bx} \frac{R^\theta}{\mH(\bx)} \int_{\Omega} |u|^p\, d\mu\, d\mH=R^\theta \|u\|^p_{L^p(\Omega, V)}.\]
For any $k\geq 0$, it follows from the doubling property of $\mu$, the upper codimension relation \eqref{upper} and the local $p$-Poincar\'e inequality  that
\begin{align}
\|T_{k+1}u- T_k u\|^p_{L^p(\bx, V)}&\leq \int_{\bx} \left(\vint_{B(z, 2^{-k-1}R)\cap \Omega} |u-u_{B(z, 2^{-k R})\cap \Omega}|\, d\mu\right)^p\, d\mH\notag\\
&\lesssim \int_{\bx} \left(\vint_{B(z, 2^{-k}R)\cap \Omega} |u-u_{B(z, 2^{-k R})\cap \Omega}|\, d\mu\right)^p\, d\mH\notag\\
&\lesssim \int_{\bx} \frac{(2^{-k}R)^p}{\mu(B(z, 2^{-k}\lambda R)\cap \Omega))}\int_{{B(z, 2^{-k}\lambda R)\cap \Omega}} g_u(x)^p\, d\mu(x)\, d\mH(z)\notag\\
&\lesssim \int_{\bx} \frac{(2^{-k}\lambda R)^{p-\theta}}{\mH(B(z, 2^{-k}R)\cap \bx))}\int_{{B(z, 2^{-k}\lambda R)\cap \Omega}} g_u(x)^p\, d\mu(x)\, d\mH(z)\notag\\
&\lesssim  \int_{\Omega(2^{-k}\lambda R)} g_u(x)^p\int_{B(x, 2^{-k}\lambda R)\cap \bx)}\frac{(2^{-k}\lambda R)^{p-\theta}}{\mH(B(z, 2^{-k}\lambda R)\cap \bx)}d\mH(z)\, d\mu(x)\label{eq-Banach-1-1}\\
&\lesssim (2^{-k} R)^{p-\theta}\int_{\Omega(2^{-k}\lambda R)} g_u(x)^p\, d\mu(x)\leq  (2^{-k}  R)^{p-\theta} \|g_u\|^p_{L^p(\Omega)},\label{tmp-1}
\end{align}
where $\Omega(r):=\{x\in \Omega: d(x, \bx)<r\}$ and the second last inequality used the fact that $\mH$ is doubling.

Since $p>\theta$, combing the estimates of $\|T_0 u\|^p_{L^{p}(\bx, V)}$ and $\|T_{k+1}u- T_k u\|^p_{L^p(\bx, V)}$, we obtain that
$$\|u^*\|_{L^p(\bx)}\lesssim R^\theta \|u\|_{L^p(\Omega, V)} + \sum_{k\geq 0} (2^{-k}  R)^{1-\theta/p}  \|g_u\|_{L^p(\Omega)}\lesssim \|u\|_{N^{1, p}(X, V)}<\infty.$$
Thus, $\widetilde T u$ exists $\mH$-almost everywhere on $\bx$. Moreover, since $|\widetilde Tu|\leq u^*$, we have 
\[\|\widetilde Tu\|_{L^p(\bx, V)}\leq \|u^*\|_{L^p(\bx)} \lesssim \|u\|_{N^{1, p}(X, V)}.\]

The proof will be complete once we show  $\widetilde Tu=Tu$ on $\bx$. For this, it suffices to show that the eqaution \eqref{trace-defn} holds with $\widetilde Tu(z)$ for $\mH$-almost every $z\in \bx$. Set 
$$E=\{z\in \bx: M_{\theta, p} g_u(z)<\infty\ \text{and}\ T_k u(z)\rightarrow \widetilde Tu(z)\ \text{as}\ k\rightarrow \infty\}.$$ 
Then Lemma \ref{fractional-maximal} implies that $\mH(E)=0$. 

For any  $0<r\leq R$, let $k_r\in \N$ such that $2^{-k_r-1}R<r\leq 2^{-k_r}R$. Then it follows from the doubling property of $\mu$ that  for any $z\in \bx\setminus E$ and $0<r\leq R$, 
\begin{align*}
\vint_{B(z, r)\cap\Omega} |u-\widetilde Tu(z)|\, d\mu&\leq \vint_{B(z, r)\cap\Omega} |u-T_{k_r}(z)|\, d\mu+|T_{k_r}(z)-\widetilde Tu(z)| \\
&\lesssim \vint_{B(z, 2^{-k_r}R)\cap \Omega}|u-u_{B(z, 2^{-k_r}R)\cap \Omega}|\, d\mu +|T_{k_r}(z)-\widetilde Tu(z)| \\
&\lesssim 2^{-k_r}R \left(\vint_{B(z, 2^{-k_r}\lambda R)\cap \Omega} g_u(x)^p\, d\mu(x)\right)^{1/p} +|T_{k_r}(z)-\widetilde Tu(z)|\\
&\lesssim (2^{-k_r}R)^{1-\theta/p} M_{\theta, p} g(z)+ |T_{k_r}(z)-\widetilde Tu(z)|.
\end{align*}
Since $z\in \bx\setminus E$ and $k_r\rightarrow \infty$ as $r\rarrow 0$, we have
$$\vint_{B(z, r)\cap\Omega} |u-\widetilde Tu(z)|\, d\mu\rarrow 0,\ \ \text{as}\ \ r\rarrow 0.$$
Hence \eqref{trace-defn} holds with $\widetilde Tu(z)$ for $\mH$-almost every $z\in \bx$. The proof is complete.
\end{proof}

\begin{proof}[Proof of Theorem \ref{thm-trace-1}]
This is a direct consequence of Theorem \ref{thm-Banach-1} and Lemma \ref{trace-isometric}.
%{\color{red} More detailed proof should be given after the metric valued Sobolev space is defined, but it should be direct since we have Lemma \ref{trace-isometric}.}
\end{proof}

%\subsection{Proof of Theorem \ref{thm-trace-1}}

As a consequence of the proof of Theorem \ref{thm-Banach-1}, we obtain the following convergence result for traces of metric valued Sobolev spaces, which in particular gives Theorem \ref{thm:convergence of trace}. 
%In the special case when $\Omega$ is a bounded Lipschitz domain in a smooth Riemannian manifold, it was proved by Korevaar and Schoen \cite[Theorem 1.12.2]{ks93}.
\begin{thm}\label{prop:convergence for trace}
	Suppose $\Omega\subset X$ is $(p,\theta)$-admissible. Let $\{u_i\}\subset N^{1,p}(\Omega,Y)$ be a sequence with uniformly bounded energy, that is, 
	$$\sup_{i\in \mathbb{N}}\int_{\Omega}g_{u_i}^pd\mu<\infty.$$
	If $u_i$ converges to some $u\in N^{1,p}(\Omega,Y)$ in $L^p(\Omega,Y)$, then $Tu_i\to Tu$ in $L^p(\bx,Y)$. Furthermore, two maps $u, v\in N^{1, p}(\Omega, Y)$ have the same trace if and only if $d(u,v)\in N^{1, p}(\Omega, \real)$ and has zero trace.
\end{thm}

%\begin{lem}
%	Let $\Omega$ be a $(p, \theta)$-admissible domain with $p>1$ and $\theta<p<\infty$. If the sequence $\{u_i\}\subset N^{1, p}(X, V)$ has uniformly bounded energies, and if $\{u_i\}$ converges to a function  $u\in N^{1, p}(X)$ in $L^p(X, V)$, then the trace function  $T u_i$ of $u_i$ converges to the trace function $Tu$ of $u$ in $L^p(\bx, V)$.
%\end{lem}
%{\color{red} In the lemma I assume that $u\in N^{1, p}$ because I need that the trace of $u$ exists. However this is not necessary once we use the reflexive property of $N^{1, p}$ and Mazur's lemma. But it seems the reflexive property is not immediately true if $V$ is just a Banach space. I am not sure what kind of lemma we want.}
\begin{proof}
	For both assertions, embedding $Y$ isometrically into some Banach space $V$ if necessary, we may assume $Y=V$ is a Banach space.
	
	For the first claim, recall that in the proof of Theorem \ref{thm-Banach-1}, we proved that $Tf= \widetilde Tf$ for any $f\in N^{1, p}(X, V)$, where
	$$\widetilde Tf =\lim_{k\rightarrow \infty}T_k f.$$
	It follows from the estimate \eqref{tmp-1} that 
	\begin{align*}
		\|Tf-T_k f\|_{L^p(\bx, V)} &\leq \sum_{j\geq k}\|T_{j+1} f-T_j f\|_{L^p(\bx, V)}\\
		& \lesssim \sum_{j\geq k} (2^{-j}R)^{1-\theta/p}\|g_f\|_{L^p(\Omega)}
		\lesssim (2^{-k} R)^{1-\theta/p} \|g_f\|_{L^p(\Omega)},
	\end{align*}
	where $g_f$ is the mimimal upper gradient of $f$.
	
	Hence for any two functions $f, h\in N^{1, p}(X, V)$ and any $k\in \N$, we have 
	\begin{align}
		\|Tf -Th\|_{L^p(\bx, V)} &\leq \|Tf-T_k f\|_{L^p(\bx, V)}+\|Th-T_k h\|_{L^p(\bx, V)}+\|T_kf-T_k h\|_{L^p(\bx, V)} \notag\\
		&\lesssim (2^{-k} R)^{1-\theta/p} \left(\|g_f\|_{L^p(\Omega)}+\|g_h\|_{L^p(\Omega)}\right)+\|T_kf-T_k h\|_{L^p(\bx, V)}, \label{tmp-2}
	\end{align}
	where $g_f$ and $g_h$ are minimal upper gradients of $f$ and $h$, respectively.
	Notice that for any $z\in \bx$, we have
	$$T_k f(z)=\vint_{B(z, 2^{-k} R)\cap\Omega}f \, d\mu\ \ \text{and}\ \ T_k h(z)=\vint_{B(z, 2^{-k} R)\cap\Omega} h\, d\mu.$$
	Thus
	\begin{align*}
		\|T_kf-T_k h\|^p_{L^p(\bx, V)} &=\int_{\bx} |T_kf(z)-T_k h(z)|^p\,d\mH(z)\\
		&\leq\int_{\bx} \left(\vint_{B(z, 2^{-k}R)\cap\Omega}|f(x)- h_{B(z, 2^{-k}R)\cap\Omega}| \, d\mu(x)\right)^pd\mH(z)\\
		&\lesssim \int_{\bx} \vint_{B(z, 2^{-k}R)\cap\Omega} |f(x)-h(x)|^p\, d\mu(x)\, d\mH(z)\\
		&\quad\ +\int_\bx  \left(\vint_{B(z, 2^{-k}R)\cap\Omega}|h(x)- h_{B(z, 2^{-k}R)\cap\Omega}| \, d\mu(x)\right)^pd\mH(z)=: I_1+I_2.
	\end{align*}
	%The estimates of $I_1$ and $I_2$ are 
	Using similar arguments as that in \eqref{tmp-1}, we obtain 
	$$I_2\lesssim (2^{-k}R)^{p-\theta}\|g_h\|^p_{L^p(\Omega)}.$$
For the estimate of $I_1$, it follows from  the upper codimension relation \eqref{upper} that
	\begin{align*}
	I_1&=\int_{\bx} \frac{1}{\mu(B(z, 2^{-k}R)\cap\Omega)} \int_{B(z, 2^{-k}R)\cap\Omega} |f(x)-h(x)|^p\, d\mu(x)\, d\mH(z)\\
	&\lesssim \int_{\bx} \frac{(2^{-k}R)^{-\theta}}{\mH(B(z, 2^{-k}R)\cap \partial\Omega)} \int_{B(z, 2^{-k}R)\cap\Omega} |f(x)-h(x)|^p\, d\mu(x)\, d\mH(z).
	\end{align*}
   Then by a similar argument with the one used in \eqref{eq-Banach-1-1} and \eqref{tmp-1}, we arrive at the estimate 
   \[I_1\lesssim (2^{-k}R)^{-\theta} \|f-h\|^p_{L^p(\Omega, V)}\]

	Thus, the estimate \eqref{tmp-2} can be rewritten as
	\begin{equation}\label{tmp-3}
		\|Tf-Th\|_{L^p(\bx, V)}\lesssim  (2^{-k} R)^{1-\theta/p} \left(\|g_f\|_{L^p(\Omega)}+\|g_h\|_{L^p(\Omega)}\right)+(2^{-k}R)^{-\theta/p} \|f-h\|_{L^p(\Omega, V)}.
	\end{equation}
	The above inequality shows that if the sequence $u_i$ converges to $u$ in $L^p(\Omega, V)$ and if the sequence has uniformly bounded energy, then  $Tu_i$ converges to $Tu$ in $L^p(\bx, V)$. Indeed, if we choose $f=u$ and $h=u_i$ in the above inequality, we know form the lower semicontinuity of energy (see \cite[Theorem 7.3.9]{hkst12}) that the energy of $u$ is also bounded and hence the first term on the right-hand side of \eqref{tmp-3} can be made arbitrary small by choosing $k$ big enough. Once $k$ is fixed, the second term can be made small by choosing $i$ large.
	
	We now turn to the second claim and assume that $u, v\in N^{1, p}(\Omega, V)$ have the same trace, i.e., $Tu(x)=Tv(x)$ for $\mH$-a.e. $x\in \bx$. We first show that $d(u,v)=|u-v|\in N^{1, p}(\Omega)$. Since $|u-v|\leq |u|+|v|$, $|u-v|\in L^p(\Omega)$. The minimal upper gradient $g_{|u-v|}$ of $|u-v|$ is controlled by $g_u+g_v$, where $g_u$ and $g_v$ are minimal upper gradients of $u$ and $v$. Indeed, for any rectifiable curve $\gamma$ connecting $x, y\in \Omega$, by triangle inequality, we have that
	\begin{align*}
		\big||u(x)-v(x)|-|u(y)-v(y)|\big|\leq |u(x)-u(y)|+|v(x)-v(y)| \leq \int_\gamma g_u+g_v\, ds.
	\end{align*}
	Thus, $|u-v|\in N^{1, p}(\Omega)$. Since $Tu(x)=Tv(x)$ for $\mH$-a.e. $x\in \bx$, it follows from the definition of trace that for $\mH$-a.e. $x\in \bx$, we have
	\begin{align*}
		\lim_{r\rightarrow 0^+} \vint_{B(x, r)\cap\Omega} |u-v|\, d\mu&\leq \lim_{r\rightarrow 0^+} \vint_{B(x, r)\cap\Omega} |u-Tu(x)|\,d\mu+ |Tu(x)-Tv(x)|\\
		&\quad\quad\quad +\lim_{r\rightarrow 0^+} \vint_{B(x, r)\cap\Omega} |Tv(x)-v(x)|\,d\mu=0.
	\end{align*}
	Hence $|u-v|$ has trace zero.
	
	For the converse, assume that $|u-v|\in N^{1, p}(\Omega)$  has trace zero. Notice that for any $x\in \bx$ and any $y\in {B(x, r)\cap\Omega}$, we have that
	$$|Tu(x)-Tv(x)| \leq |Tu(x)-u(y)|+|u(y)-v(y)|+|v(y)-Tv(x)|.$$
	It follows from the definition of trace that for $\mH$-a.e. $x\in \bx$, we have
	\begin{align*}
		|Tu(x)-Tv(x)|&\leq \vint_ {B(x, r)\cap\Omega} |Tu(x)-u|\, d\mu+ \vint_ {B(x, r)\cap\Omega} |u-v|\, d\mu\\
		&\quad\quad\quad+ \vint_ {B(x, r)\cap\Omega} |v-Tv(x)|\, d\mu\rightarrow 0, \ \text{as}\ r\rightarrow 0.
	\end{align*}
	Thus, $u$ and $v$ have the same trace.
\end{proof}

%\begin{lem}
%	Two functions $u, v\in N^{1, p}(\Omega, V)$ have the same trace if and only if $|u-v|\in N^{1, p}(\Omega, \real)$  has trace zero.
%\end{lem}
%\begin{proof}
%
%\end{proof}

%{\color{blue}
%We additionally require that 
%
%\begin{itemize}
%	\item For each $u\in N^{1,p}(\Omega)$ that can be approximated in $\|\cdot\|_{N^{1,p}}$ by $N^{1,p}_c(\Omega)$-functions, the zero extension belongs to $N^{1,p}(X)$;
%	
%	\item For each $u\in N^{1,p}(\Omega)$ with zero trace, the zero extension belongs to  $N^{1,p}(X)$. 
%\end{itemize}
%}
%
%As a consequence of the above requirement, we have
%
%\begin{lem}\label{lemma:Poincare zero boundary value}
%	If $u\in N^{1,p}(\Omega)$ with $Tu=0$ $\mH$-almost everywhere on $\bx$, then 
%	\begin{equation}\label{eq:Poincare for Sobolev functions with zero boundary value 2}
%		\|u\|_{L^{p}(\Omega)}\leq C(\Omega)\|g_u\|_{L^p(\Omega)}.
%	\end{equation}
%\end{lem}
%
%Discuss $N^{1,p}=M^{1,p}$ and extension in more details.

\section{Solution to the Dirichlet problem}\label{sec:solution to Dirichlet}

\subsection{Hajlasz-Sobolev spaces and consequences}

Let $\Omega\subset X$ be a domain and $Y$ a complete metric space.
\begin{defn}[Hajlasz-Sobolev spaces]\label{def:Hajlasz-Sobolev}
	A measurable map $u\colon \Omega\to Y$ belongs to the Hajlasz-Sobolev space $M^{1,p}(\Omega,Y)$ if $u\in L^p(X,Y)$ and there exists a nonnegative function $g\in L^p(\Omega)$ such that the Hajlasz gradient inequality
	\begin{equation}\label{eq:Hajlasz gradient}
		d_Y(u(x),u(z))\leq d_X(x,z)\left(g(x)+g(z)\right)
	\end{equation} 
   holds for all $x,y\in \Omega\backslash N$ for some $N\subset \Omega$ with $\mu(N)=0$. For each $u\in M^{1,p}(\Omega,Y)$, the associated Hajlasz energy $E^p_H(u)$ is defined as 
   $$E^p_H(u):=\inf_{g}\|g\|_{L^p(\Omega)},$$
   where the infimum is taken over all Hajlasz gradient $g$ of $u$, that is, $g$ such that \eqref{eq:Hajlasz gradient} holds.
\end{defn}

The following equivalence of metric valued Sobolev spaces is well-known.
\begin{prop}[\cite{hkst12}, Corollary 10.2.9]\label{prop:equivalence of M and N}
 Suppose $\mu$ is doubling and $\Omega$ supports a $q$-Poincar\'e inequality for some $1\leq q<p$. Then $M^{1,p}(\Omega,Y)=N^{1,p}(\Omega,Y)$. Furthermore, there exists a constant $C\geq 1$, depending only on the data associated to $\Omega$, such that for each $u\in M^{1,p}(\Omega,Y)=N^{1,p}(\Omega,Y)$, 
 $$C^{-1}E^p(u)\leq E^p_H(u)\leq CE^p(u).$$
\end{prop}

%\begin{rem}
%	content...
%\end{rem}

The proof of Theorem \ref{thm:ultra-completion} requires the following Rellich compactness result, which was proved in~\cite[Theorem 3.1]{gw17} when $\Omega\subset \R^n$ is a bounded Lipschitz domain.
\begin{thm}[Generalized Rellich compactness]\label{thm:Rellich compactness}
Suppose $\mu$ is doubling and $\Omega$ supports a $q$-Poincar\'e inequality for some $1\leq q<p$.	For every $m\in \mathbb{N}$, let $(Y_m,d_m)$ be a complete metric space, $K_m\subset Y_m$ compact and $\{u_m\}\subset N^{1,p}(\Omega,Y_m)$. Suppose that $(K_m,d_m)$ is uniformly compact and
	\begin{align}\label{eq:generalized compactness}
		\sup_{m\in \mathbb{N}}\int_{\Omega}d_m^p(u_m(x),y_m)d\mu(x)+E^p(u_m)<\infty 
	\end{align}
	for some and thus every $y_m\in K_m$. Then after possibly passing to a subsequence, there exist a complete metric space $Z$, a compact subset $K\subset Z$, isometric embeddings $\varphi_k\colon Y_m\to Z$, and $v\in N^{1,p}(\Omega,Z)$ such that $\varphi_m(K_m)\subset K$ for all $m\in \mathbb{N}$ and $\varphi_m\circ u_m$ converges to $v$ in $L^p(\Omega,Z)$.
\end{thm}

Recall that a sequence of compact metric space $(B_m,d_m)$ is called uniformly compact if $\sup_{m}\diam B_m<\infty$ and if for every $\varepsilon>0$, there exists $N\in \mathbb{N}$ such that every $B_m$ can be covered by at most $N$ balls of radius $\varepsilon$.

\begin{proof}[Proof of Theorem \ref{thm:Rellich compactness}]
	The proof is essentially contained in~\cite[Theorem 3.1]{gw17} and thus we only point out the necessary changes.  The key ingredient of the proof is that under our assumptions on $\Omega$, we have by Proposition \ref{prop:equivalence of M and N} that $N^{1,p}(\Omega,Y)= M^{1,p}(\Omega,Y)$ and each $u_m\in N^{1,p}(\Omega,Y)$ satisfies the pointwise inequality \eqref{eq:Hajlasz gradient}
	\begin{equation*}\label{eq:hajlasz estimate}
		d_m(u_m(x),u_m(x'))\leq d(x,x')(h_m(x)+h_m(x'))
	\end{equation*}
	almost everywhere for some $h_m\in L^p(\Omega)$ with $\|h_m\|_{L^p(\Omega)}^p\leq C(p,\Omega)E^p(u_m)$. 
\end{proof}

\subsection{Proof of Theorem \ref{thm:main existence} }

In this section, we provide the proof of Theorem \ref{thm:main existence}, which is very similar to \cite[Proof of Theorem 1.4]{gw17}. In the first step, we prove the following result on ultra-limits of subsequences of Sobolev maps, which extends \cite[Theorem 1.6]{gw17}.

\begin{thm}\label{thm:ultra-completion}
Suppose $\Omega\subset X$ is a $(p,\theta)$-admissible domain and $Y_\omega$ is an ultra-completion of the complete metric space $Y$. If $\{u_k\}\subset N^{1,p}(\Omega,Y)$ is a bounded sequence for some $p>1$, then, after possibly passing to a subsequence, the map $\phi(z):=[(u_m(z))]$ belongs to $W^{1,p}(\Omega,Y_\omega)$ and satisfies
$$E^p(\phi)\leq \liminf_{k\to \infty}E^p(u_k).$$
Moreover, if $Tu_k$ converges to some  map $\rho\in L^p(\bx,Y)$ $\mH$-almost everywhere on $\bx$, then $T\phi=\iota\circ \rho$. 
\end{thm}

\begin{proof}[Proof of Theorem~\ref{thm:ultra-completion}]
 The proof is essentially contained in \cite[Proof of Theorem 1.6]{gw17} and we present it again for the convenience of the readers. After possibly passing to a subsequence, we may assume that
 $$E^p(u_k)\to\liminf_{m\to\infty}E^p(u_m)$$
 as $k\to \infty$.
 
 Fix $y_0\in Y$ and apply the Rellich compactness Theorem \ref{thm:Rellich compactness}. After possibly passing to a subsequence, there exist a complete metric space $Z=(Z,d_Z)$, a compact subset $K\subset Z$, and isometric embeddings $\varphi_k\colon Y\to Z$ and $v\in N^{1,p}(\Omega,Z)=M^{1,p}(\Omega,Z)$ such that $\varphi_k(y_0)\subset K$ for all $k$ and
 %$\varphi_k(C_m)\subset Z^m$ for all $k,m$ and 
 $v_k:=\varphi_k\circ u_k$ converges in $L^p(\Omega,Z)$ to $v$ as $k\to \infty$. 
 After passing to a further subsequence, we may assume that $v_k$ converges almost everywhere to $v$ on $\Omega$. Let $N\subset \Omega$ be a set of $\mu$-measure zero such that $v_k(z)\to v(z)$ for all $z\in \Omega\backslash N$.

   Define a subset of $Z$ by $B:=\{v(z):z\in \Omega\backslash N \}$.
	%$\cup \varphi(C)$. 
	The map $\psi\colon B\to Y_\omega$, given by $\psi(v(z))=[(u_k(z))]$ when $z\in \Omega\backslash N$
	% and $\psi(\varphi(y))=[(y)]$ when $y\in C$ 
	is well-defined and isometric by~\cite[Lemma 2.2]{gw17}. Since $Y_\omega$ is complete, there exists a unique extension of $\psi$ to $\overline{B}$, which we denote again by $\psi$. After possibly redefining the map $v$ on $N$, we may assume that $v$ has image in $\overline{B}$ and hence $v$ is an element of $M^{1,p}(\Omega,\overline{B})$. Now, we define a mapping by $$\phi(z):=\psi(v(z))=[(u_k(z))]$$ 
	and then $\phi$ belongs to $M^{1,p}(\Omega,Y_\omega)$ and by the lower semicontinuity of upper gradient energy \cite[Theorem 7.3.9]{hkst12} it satisfies
	\begin{align}\label{eq:5}
		E^p(\phi)\leq E^p(v)\leq \liminf_{k\to \infty}E^p(v_k)=\lim_{k\to \infty} E^p(u_k).
	\end{align}
	
	It remains to prove the trace equality. Suppose $Tu_k$ converges to some map $\rho\in L^p(\bx,Y)$ almost everywhere on $\bx$. Arguing as in \cite[Page 104]{gw17}, we can find compact subsets $C_1\subset C_2\subset \cdots \subset Y$, isometric embeddings $\varphi_k\colon Y\to Z$ and $v\in N^{1,p}(\Omega,Z)$ such that $v_k:=\varphi_k\circ u_k$ converges in $L^p(\Omega,Z)$ to $v$ as $k\to \infty$. Furthermore, if we set $C=\bigcup_{l=1}^\infty C_l$, then aftering passing to a further subsequence if necessary we may assume that $v_k$ converges to $v$ almost everywhere on $\Omega$ and $\varphi_k|_{C}$ converges pointwise to an isometric embedding $\varphi\colon C\to Z$, with the convergence being uniform on each $C_k$. Let $N\subset \Omega$ be a set of $\mu$-measure zero such that $v_k(z)\to v(z)$ for all $z\in \Omega\backslash N$. 
	
	Define a subset of $Z$ by
	$$B:=\{v(z):z\in\Omega \}\cup \varphi(C).$$
    The map $\psi\colon B\to Y_\omega$ given by 
    \begin{equation*}
    	\begin{cases}
    	\psi(v(z))=[(u_k(z))] & \text{ if } z\in \Omega\backslash N,\\
    	\psi(\varphi(x))=\iota(x)=[(x)] & \text{ if }x\in C,
    	\end{cases}
    \end{equation*}
    is well-defined and an isometric embedding by \cite[Lemma 2.2]{gw17}. Since $Y_\omega$ is complete, there exists a unique isometric extension of $\psi$ to $\overline{B}$, which we denote again by $\psi$. After possibly redefining the map $v$ on $N$, we may assume $v\in M^{1,p}(\Omega,\overline{B})$. The map $\phi(z):=\psi(v(z))=[(u_k(z))]$ then belongs to $M^{1,p}(\Omega,Y_\omega)$ and satisfies \eqref{eq:5}. Moreover, by Proposition \ref{prop:convergence for trace}, we have that $Tv_k=\varphi_k\circ Tu_k$ converges to $\varphi\circ \rho$ almost everywhere on $\bx$ and a subsequence of $Tv_k$ converges to $Tv$ almost everywhere. It thus follows that $Tv=\varphi\circ \rho$ and hence
    $$T\phi=\psi\circ Tv=\psi\circ\varphi\circ \rho=\iota\circ \rho.$$ 
    The proof is complete.
    
\end{proof}

With Theorem \ref{thm:ultra-completion} at hand, the proof of Theorem \ref{thm:main existence} is immediate. 
\begin{proof}[Proof of Theorem \ref{thm:main existence}]
	Let $\phi\in N^{1,p}(\Omega,Y)$ and let $\{u_k\}\subset N^{1,p}(\Omega,Y)$ be an energy minimizing sequence with $Tu_k=T\phi$ for each $k$. Then by the characterization of trace from Proposition \ref{prop:convergence for trace}, $h_k(x)=d(u_k(x),\phi(x))\in N^{1,p}_0(\Omega)$.  Since $\sup_{k}E^p(h_k)<\infty$, it follows from the global $p$-Poincar\'e inequality \eqref{eq:Poincare for Sobolev functions with zero boundary value} that $\sup_{k}\|h_k\|_{L^p(\Omega)}<\infty$. Hence
	\begin{align*}
		\sup_{k}\int_{\Omega} d^p(u_k(x),y_0)d\mu(x)+E^p(u_k)<\infty.
	\end{align*}
	Thus $\{u_k\}$ is a bounded sequence in $N^{1,p}(\Omega,Y)$. Let $Y_\omega$ be an ultra-completion of $Y$ such that $Y$ admits a 1-Lipschitz retraction $P\colon Y_\omega\to Y$. After possibly passing to a subsequence, we may assume by Theorem \ref{thm:ultra-completion} that the map $v(z):=[(u_k(z))]$ belongs to $N^{1,p}(\Omega,Y_\omega)$ and satisfies $Tv=\iota\circ T\phi$ and 
	$$E^p(v)\leq \lim_{k\to \infty}E^p(u_k).$$
	Since $P\colon Y_\omega\to Y$ is a 1-Lipschitz retraction, the map $u:=P\circ v$ belongs to $N^{1,p}(\Omega,Y)$ and satisfies $Tu=T\phi$ and $E^p(u)\leq \lim_{k\to\infty}E^p(u_m)$. The proof is complete.	
\end{proof}

\section{The theory of trace in the borderline case}\label{sec:refined trace}

In Theorem \ref{thm-Banach-1}, $\Omega$ is assumed to be weakly $(p,\theta)$-admissible, where in Definition \ref{def:p admissible domain} it is required that $p>\theta$. It is natural to ask for what happens if $p=\theta$. We shall address this problem in this section. To simplify our exposition, we say that $\Omega$ is weakly $\theta$-admissible if it is weakly $(p,\theta)$-admissible with $p=\theta$. 

For the borderline case when $\Omega$ is weakly $\theta$-admissible for some $\theta>1$, the traces of $N^{1, p}(X, V, d\mu)$ may not exist. To characterize the existence of the traces, we give an addition weight  on the measure $\mu$ and investigate the relationship between the existence of the traces and the properties of the weight function.

For a locally integrable weight function $\rho\colon \Omega\rightarrow (0, \infty)$, we  define $L^p(\Omega, V, \rho d\mu)$ and $N^{1, p}(\Omega, V, \rho d\mu)$ by replacing  $d\mu$ with the weighted measure $\rho \mu$ in the integrals of the norms. 

The following result follows essentially from the proof of Theorem \ref{thm-Banach-1}, whose core idea is  similar with the one used in  \cite{m17}
\begin{thm}\label{thm-Banach-2}
Suppose that $\Omega$  is $\theta$-admissible for some $\theta>1$. Let $\omega\colon (0, \infty)\rightarrow[1, \infty)$ be a non-increasing function and  $\rho(x):=w(\dist(x, \bx))$. Then the trace operator 
$$T\colon N^{1, \theta}(\Omega, V, \rho d\mu)\rightarrow L^\theta(\partial \Omega, V)$$ 
is bounded and linear, provided that $\int_0^1 {t^{-1} \omega(t)^{-1/\theta}}\, dt<\infty$.
\end{thm}
\begin{proof}
The proof is a minor modification of that used in Theorem \ref{thm-Banach-1} and thus we use the same notations here.
Repeat the proof of Theorem \ref{thm-Banach-1} until the estimate \eqref{eq-Banach-1-1}. Since $\omega$ is non-increasing, by considering the weighted measure $\rho(x)d\mu$ instead of the measure $d\mu$ in the  the estimate \eqref{eq-Banach-1-1}, we obtain that
\begin{align*}
\|T_{k+1}u- &T_k u\|^\theta_{L^\theta(\bx, V)}\\
&\lesssim \int_{\Omega(2^{-k}\lambda R)} g(x)^\theta\rho(x)\int_{B(x, 2^{-k} \lambda R)\cap \bx)}\frac{d\mH(z)}{\omega(2^{-k}\lambda R)\mH(B(z, 2^{-k}\lambda R)\cap \bx)}\, d\mu(x)\\
&\lesssim \frac{1}{\omega(2^{-k}\lambda R)} \int_{\Omega(2^{-k}\lambda R)} g(x)^\theta\rho(x)\, d\mu=\frac{1}{\omega(2^{-k}\lambda R)} \|g\|^\theta_{L^\theta(\Omega, \rho d\mu)}.
\end{align*}
Since $\omega(t)\geq 1$, we have 
$$\|T_0 u\|^\theta_{L^\theta(\Omega)}\lesssim R^\theta\|u\|^\theta_{L^\theta(\Omega, V)}\leq R^\theta\|u\|^\theta_{L^\theta(\Omega, V, \rho d\mu)}.$$
Thus, it follows from $\int_0^1 t^{-1}w(t)^{-1/\theta}\, dt<\infty$ that 
\begin{align*}
\|u^*\|_{L^\theta(\Omega)}&\lesssim R^\theta \|u\|_{L^\theta(\Omega, V, \rho d\mu)} + \sum_{k\geq 0} \frac{1}{(\omega(2^{-k}\lambda R))^{1/p}}  \|g\|_{L^\theta(\Omega, \rho d\mu)}\\
&\approx  \|u\|_{L^\theta(\Omega, V, \rho d\mu)} +\sum_{k\geq 0} \frac{\|g\|_{L^\theta(\Omega, \rho d\mu)}}{(\omega(2^{-k}\lambda R))^{1/\theta}}  \int_{2^{-k}\lambda R}^{2^{-k+1}\lambda R} \frac{dt}{t}\\
&\leq \|u\|_{L^\theta(\Omega, V, \rho d\mu)} +\|g\|_{L^\theta(\Omega, \rho d\mu)} \int_{0}^{2\lambda R} \frac{dt}{t w(t)^{1/\theta}}\\
&\approx \|u\|_{L^\theta(\Omega, V, \rho d\mu)} +\|g\|_{L^\theta(\Omega, \rho d\mu)} \int_{0}^{1} \frac{dt}{t w(t)^{1/\theta}}\\
&\lesssim \|u\|_{N^{1, \theta}(X, V, \rho d\mu)}.
\end{align*}
Hence $\widetilde T u$ exists $\mH$-almost everywhere on $\bx$ and we have the estimate 
\[\|\widetilde Tu\|_{L^\theta(\bx, V)}\leq \|u^*\|_{L^\theta(\bx)} \lesssim \|u\|_{N^{1, \theta}(X, V, \rho d\mu)}.\]

The proof will be complete once we show  $\widetilde Tu=Tu$. We shall use the similar idea as in the proof of Theorem \ref{thm-Banach-1}, and set
$$E=\{z\in \bx: M_{\theta, \theta} (g\rho^{1/\theta})(z)<\infty\ \text{and}\ T_k u(z)\rightarrow \widetilde Tu(z)\ \text{as}\ k\rightarrow \infty\}.$$
Then $\mH(E)=0$ and for every $z\in \bx\setminus E$ and every $2^{-k_r-1}R<r\leq 2^{-k_r}R$, we have
\begin{align*}
\vint_{B(z, r)\cap\Omega} |u-\widetilde Tu(z)|\, d\mu&\lesssim 2^{-k_r}R \vint_{B(z, 2^{-k_r}\lambda R)\cap \Omega} g(x)^p\, d\mu(x) +|T_{k_r}-\widetilde Tu(z)|\\
&\leq \frac{2^{-k_r}R}{\omega(2^{-k_r}\lambda R)}\left(\vint_{B(z, 2^{-k_r}\lambda R)\cap \Omega} g(x)^\theta \rho(x)\, d\mu(x)\right) +|T_{k_r}-\widetilde Tu(z)|\\
&\leq \frac{1}{w(\lambda r/2)} M_{\theta, \theta} (g\rho^{1/\theta})(z)+|T_{k_r}-\widetilde Tu(z)|.
\end{align*}
Since $\int_0^1 t^{-1}w(t)^{-1/\theta}\, dt<\infty$, the non-increasing function $w(t)\rightarrow \infty$ as $t\rightarrow 0$, and so it follows 
$$\vint_{B(z, r)\cap\Omega} |u-\widetilde Tu(z)|\, d\mu\rarrow 0,\ \ \text{as}\ \ r\rarrow 0.$$
This means that \eqref{trace-defn} holds with $\widetilde Tu(z)$ for $\mH$-almost every $z\in \bx$ and the proof is thus complete.
\end{proof}

The integrability condition on $\omega$ can be relaxed if $\Omega$ has certain nice geometry. For this, we shall introduce the so-called John domains, which plays an important role in geometric analysis in metric spaces \cite{hak00}.

\begin{defn}
A bounded domain $\Omega\subset X$ is called a {\it John domain} with {\it John constant} $c_J\in (0, 1]$ and {\it John center} $a\in \Omega$ if every $x\in \Omega$ can be joined to $a$ by a rectifiable curve $\gamma: [0, \ell_\gamma]\rightarrow \Omega$ parametrized by arc-length such that $\gamma(0)=x$, $\gamma(\ell_\gamma)=a$ and
\begin{equation}\label{John-condition}
\dist(\gamma(t), X\setminus\Omega)\geq c_J t\ \ \ \text{for all }\ t\in [0, \ell_\gamma].
\end{equation}
\end{defn}

If $\Omega$ is a John domain with compact closure in $X$, the it follows from the Arzela-Ascoli theorem that every $z\in \bx$ can be joined to the John center by a rectifiable curve such that \eqref{John-condition} holds.
\begin{thm}\label{thm-Banach-John}
Let $\Omega$,  $\theta$, $w$ and $\rho$ be as in Theorem \ref{thm-Banach-2}. 
Assume additionally that $\Omega$ is a John domain with compact closure. Then the trace operator $$T\colon N^{1, \theta}(\Omega, V, \rho d\mu)\rightarrow L^\theta(\partial \Omega, V)$$ 
is bounded and linear provided
\begin{equation}\label{eq:John}
\int_0^1 {t^{-1} \omega(t)^{-1/(\theta-1)}}\, dt<\infty.
\end{equation}
\end{thm}

\begin{proof}
We first show that the trace operator $T$ is bounded and linear if \eqref{eq:John} holds. Let $u\in N^{1, \theta}(\Omega, V)$ be fixed.
 Let $\delta=\dist(a, \bx)>0$, where $a\in \Omega$ is the John center of $\Omega$. For any point $z\in \bx$, let $\gamma_z$ be the arc-length parametrized curve that connects the points $a$ and $z$, and satisfies \eqref{John-condition}.

Let $t_k=\delta(1-\frac{c_J}{2\lambda})^k$ and $r_k=\frac{c_J}{2\lambda}t_k$ for $k=0, 1, \cdots$. Next, we define a chain of balls $\{B^k_z\}_{z\in\mathbb N}$ by setting
$$B_z^k=B(\gamma_z(t_k), r_k),\ \ \ \text{for all}\ \ k=0, 1,\cdots.$$
The chain $\{B^k_z\}_{z\in\mathbb N}$ consists of balls of bounded overlap, where the upper bound on the number of overlapping balls depends only on the John constant $c_J$. Moreover, abbreviating $\alpha:=(2-\frac{c_J}{2\lambda})$, we have $B_z^{k+1}\subset \alpha B_z^k$. By triangle's inequality and direct computations, for any $x\in \alpha\lambda B_z^k$, we have
\begin{equation}\label{e1}
\beta_1 r_k=:\frac{c_J r_k}{2}\leq d(x, z)\leq (\alpha\lambda+\frac{2\lambda}{c_J})r_k:=\beta_2 r_k
\end{equation}
and
\begin{equation}\label{e2}
\beta_1 r_k=\frac{c_J r_k}{2}\leq d(x, \partial\Omega)\leq (\alpha\lambda+\frac{2\lambda}{c_J})r_k=\beta_2 r_k.
\end{equation}
Then we define $T_k u(z):=\vint_{B_z^k} u\, d\mu$. 

In the next step, we will show that $\widetilde T u(z):=\lim_{k\rarrow+\infty} T_k u(z)$ exists for $\mH$-almost every $z\in \bx$ and that $T u=\widetilde T u\in L^p(\bx, V)$. It suffices to show that the function 
$$u^*=\sum_{k\geq 0} |T_{k+1}u -T_k u|+|T_0 u|=: |T^* u|+|T_0 u|$$
belongs to $L^p(\bx)$. We claim that 
$$\|u^*\|_{L^p(\bx)}\leq \|T_0 u\|_{L^p(\bx, V)}+\|T^* u\|_{L^p(\bx)}<\infty.$$
By the doubling property of $\mu$ and the relation \eqref{e1}, we know that
\begin{eqnarray*}
\|T_0 u\|_{L^p(\bx)}^\theta&\leq&\int_{\bx} \vint_{B_z^0}|u(x)|^\theta\, d\mu(x)\, d\mH(z)=\int_{\bx}\int_{\Omega} \frac{|u(x)|^\theta \chi_{B_z^0}}{\mu(B_z^0)}\, d\mu(x)\, d\mH(z)\\
&\approx&\int_{\bx} \int_{\Omega}\frac{  |u(x)|^\theta \chi_{B_z^0}}{\mu(B(z, \beta_2\delta)\cap\Omega)} \, d\mu(x)\, d\mH(z)\\
&\leq& \int_{\Omega} |u(x)|^\theta \int_{\bx} \frac{ \chi_{B(x, \beta_2\delta)\cap\bx}}{\mu(B(z, \beta_2\delta)\cap\Omega)} \, d\mH(z)\, d\mu(x)\\
&\approx& \int_{\Omega} |u(x)|^\theta \int_{B(x, \beta_2\delta)\cap\bx}\frac{\delta^{-\theta}\,d\mH(z)}{\mH(B(z, \beta_2\delta)\cap\Omega)}\, d\mu(x)\\
&\lesssim& \delta^{-\theta} \int_{\Omega} |u(x)|^\theta\, d\mu(x) \leq \delta^{-\theta} \int_{\Omega} |u(x)|^\theta \rho(x)\, d\mu(x)\approx \|u\|^\theta_{L^\theta(\Omega), V, \rho d\mu)}.
\end{eqnarray*}
Here we used the fact that  $\mH$ is doubling and $B(x, \beta_2\delta)\cap\Omega\subset B(z, 2\beta_2\delta)\cap\Omega$ whenever $z\in B(x, \beta_2\delta)$.

We next estimate the $L^p$-norm of $|T^*u|$. For any $k\geq 0$,  by the local $\theta$-Poincar\'e inequality and the fact that $B_z^{k+1}\subset B_z^k$, we know that
\begin{eqnarray*}
&&|T_k u(z)-T_{k+1} u(z)|=|u_{B_z^k}-u_{B_z^{k+1}}|\leq |u_{B_z^k}-u_{\alpha B_z^k}|+|u_{\alpha B_z^k}-u_{B_z^{k+1}}|\\
&&\ \ \ \ \ \ \ \ \ \ \ \lesssim r_k \left(\vint_{\alpha\lambda  B_z^k} g(x)^\theta\, d\mu(x)\right)^{1/\theta} 
\lesssim r_k \omega(\beta_2r_k)^{-1/\theta} \left(\vint_{\alpha\lambda  B_z^k} g(x)^\theta \rho(x)\, d\mu(x)\right)^{1/\theta}.
\end{eqnarray*}
Then, by H\"older's inequality, we have
\begin{eqnarray*}
\left(\sum_{k=0}^{\infty}|T_k u(z)-T_{k+1} u(z)|\right)^\theta \lesssim \left(\sum_{k=0}^{\infty} \omega(\beta_2 r_k)^{-1/(\theta-1)}\right)^{\theta-1}\left(\sum_{k=0}^{\infty} {r_k}^\theta \vint_{\alpha\lambda  B_z^k} g(x)^\theta \rho(x)\, d\mu(x)\right).
\end{eqnarray*}
Since $\omega$ satisfies 
$$\sum_{k=0}^{\infty} \omega(\beta_2 r_k)^{-1/(\theta-1)} \approx C+\sum_{k=1}^{\infty} \omega(\beta_2 r_k)^{-1/(\theta-1)}\int_{\beta_2 r_k}^{\beta_2 r_{k-1}}\frac{dt}{ t} \lesssim C+\int_{0}^{1} \frac{dt}{t \omega(t)^{1/(\theta-1)}}<+\infty,$$
it follows from the relation \eqref{e2} that
\begin{eqnarray*}
\|T^* u\|_{L^\theta(\bx)}^\theta&=&\int_{\bx}\left(\sum_{k=0}^{\infty}|T_k u(z)-T_{k+1} u(z)|\right)^\theta\, d\mH(z)\\
&\lesssim& \int_{\bx} \sum_{k=1}^{\infty}  {r_k}^\theta \vint_{\alpha\lambda  B_z^k} g(x)^\theta \rho(x)\, d\mu(x)\, d\mH(z)\\
&\leq& \sum_{k=0}^{\infty}\int_{\Omega({\beta_2 r_k})\setminus \Omega({\beta_1 r_k})}g(x)^\theta \rho(x) \int_{\bx}\frac{{r_k}^\theta\chi_{B(x, \beta_2 r_k)\cap\bx}}{\mu(B(z, \beta_2 r_k)\cap\Omega)} d\mH(z)\, d\mu(x)\\
&\lesssim& \sum_{k=0}^{\infty}\int_{\Omega({\beta_2 r_k})\setminus \Omega({\beta_1 r_k})}g(x)^\theta \rho(x)\int_{B(x, \beta_2 r_k)\cap\bx} \frac{ d\mH(z)}{\mH(B(z, \beta_2 r_k)\cap\bx)}\, d\mu(x)\\
&\lesssim& \sum_{k=0}^{\infty}\int_{\Omega({\beta_2 r_k})\setminus \Omega({\beta_1 r_k})}g(x)^\theta \rho(x)\, d\mu(x)\\
&\lesssim& \int_{\Omega} g(x)^\theta \rho(x) \, d\mu(x) = \|g\|^p_{L^p(\Omega, \rho d\mu)}.
\end{eqnarray*}
Hence $\widetilde T u$ exist $\mH$-almost everywhere on $\bx$ and 
\[\|\widetilde Tu\|_{L^\theta(\bx, V)}\leq \|u^*\|_{L^\theta(\bx)} \lesssim \|u\|_{N^{1, \theta}(X, V, \rho d\mu)}.\]

It is left to show  $\widetilde Tu=Tu$. This follows by a similar argument as in the proof of Theorem \ref{thm-Banach-2}, upon noticing that the convergence assumption $\int_0^1 t^{-1}w(t)^{-1/(\theta-1)}\, dt<\infty$ implies the non-increasing function $w(t)\rightarrow \infty$ as $t\rightarrow 0$. The proof is thus complete.
%%
%%Let $\tilde g(x)=g(x) w_\epsilon^{1/p}(x)$. Then $\tilde g(x)\in L^p(\Omega)$.
%%Let $E=\{z\in\bx: M_{\theta, p} {\tilde g}(z)<\infty \ \ \text{and} \ \ T_k u\rightarrow T u \ \ \text{as} \ \ k \rarrow \infty\}$. Then $\mH(\bx\setminus E)=0$. For every $z\in E$ and $r>0$, denote by $k_r$ the constant satisfying $\beta_2 r_{k_r}\leq r<\beta_2 r_{k_r-1}$. Then we obtain that
%%\begin{eqnarray*}
%%&&\vint_{B(z, r)} |u(x)-T u(x)|\, d\mu(x) \leq \vint_{B(z, r)} |u(x)-T_{k_r} u(z)|\, d\mu(x) +|T_{k_r} u(z)- T u(z)|\\
%%&&\ \ \ \ \ \   \ \leq \vint_{B(z, r)} |u(x)-u_{B(z, r)}\, d\mu(x)| +|u_{B(z, r)}-T_{k_r}u(z)| +|T_{k_r} u(z)- T u(z)|\\
%%&&\ \ \ \ \ \  \ \lesssim r\left(\vint_{B(z, r)} g^p \, d\mu\right)^{1/p}+|T_{k_r} u(z)- T u(z)|\leq \frac{M_{\theta, p}{\tilde g(z)}}{w(r)^{1/p}} +|T_{k_r} u(z)- T u(z)|,
%\end{eqnarray*}
%which approaches $0$ as $r\rarrow 0$. Thus the trace operator $T u$ satisfies (1.2).
\end{proof}

We recall the following well-known lemma from \cite{V85}.
\begin{lem}[\cite{V85}]\label{lemma3.1}
Let $(K, d_K, \mu_K)$ be a $\sigma$-finite metric measure space. Then the following conditions on $(K, d_K, \mu_K)$ are equivalent:

\item (i) $L^p(K)\subset L^q(K)$ for all $p, q\in (0, \infty)$ with $p>q$;

\item (ii) $\mu_K(K)<+\infty$.
\end{lem}

The integrability assumption \eqref{eq:John} is sharp as the following example demonstrates.
\begin{example}\label{exam:sharpness}\rm
	Let  $\Omega=\mathbb{D}\in \real^2$  with the measure $\mu$ in $\Omega$ given by $d\mu=\dist(x, \partial\mathbb D)^{\theta-1}$, $\theta>1$. On the boundary $\partial \Omega$, let $\mH$ be the $1$-dimensional Hausdorff measure. By \cite[Theorem 3.4]{Hr90}, $\Omega$ supports a local $1$-Poincar\'e inequality and hence supports a local $\theta$-Poincar\'e inequality. The doubling property of $\mu$ and upper codimension-$\theta$ regularity of $\mH$ follow by direct computations. Thus $\Omega$ is weakly $\theta$-admissible.  Let $w$ and $\rho$ be as in Theorem \ref{thm-Banach-2} with
	\begin{equation}\label{infty}
\int_0^1 {t^{-1} \omega(t)^{-1/(\theta-1)}}\, dt=\infty.
\end{equation}
Then there exists a function $u\in N^{1, \theta}(\Omega, \rho d\mu)$ such that $Tu(z)$ does not exist for all $z\in \partial\Omega$. 

\begin{proof}
To find such a function $u$, it suffices to construct a function $f\colon (0, 1)\rightarrow [0, \infty)$ such that
\begin{equation}\label{eq-example}\left\{
\begin{array}l
\int_{0}^{1} f(t) dt=+\infty,\\
\int_{0}^{1} f(t)^\theta t^{\theta-1}\omega(t)\, dt<+\infty.
\end{array}\right.
\end{equation}
Indeed, if such a function $f$ exists, then we may define $u$ by setting $u(0)=0$ and 
\begin{equation}{\label{eq-example-u}}
u(x)=\int_0^{|x|} f(t)\, dt.
\end{equation}
Then the Borel function $g\colon \Omega\rightarrow [0, \infty)$ given by $g(x)=f(|x|)$ is an upper gradient of $u$. The relation \eqref{eq-example} implies that
\begin{align*}
\|g\|^\theta_{L^\theta(\Omega, \rho d\mu)} &=\int_{\Omega} g(x)^\theta \dist(x, \bx)^{\theta-1}\omega(\dist(x, \bx))\, dx\\
&\lesssim \int_{0}^1 f(t)^\theta t^{\theta-1} w(t)\, dt<\infty.
\end{align*}
Thus it follows from the H\"older inequality that for any $t_0\in (0, 1)$, 
\[\int_{t_0}^1 f(t)\, dt\leq \left(\int_{t_0}^1 f(t)^\theta t^{\theta-1} \omega(t)\, dt\right)^{1/\theta}\left(\int_{t_0}^1 t^{-1} w(t)^{-1/(\theta-1)}\, dt\right)^{(\theta-1)/\theta}<\infty,\]
since $w(t)\geq 1$ implies that $\int_{t_0}^1 t^{-1} w(t)^{-1/(\theta-1)}\, dt\leq \int_{t_0}^1 t^{-1}\, dt<\infty$. Hence for any $t_0\in(0, 1)$, we have 
\begin{equation}
\int_{0}^{t_0} f(t)\, dt =\infty.
\end{equation}
Then for the function $u$ defined by \eqref{eq-example-u}, the trace  $Tu(z)$ does not exist for all $z\in \partial\Omega$, since as $x$ goes to the boundary $\bx$, the function $u(x)$ goes to infinity uniformly. 

As one might notice, there is a gap here, that is we do not know if the function $u$ belongs to $L^\theta(\Omega, \rho d\mu)$ or not. But this gap could be fixed by modifying the function to be an oscillatory function  with values in $[0,1]$, instead of an increasing function with respect to $|x|$, such that there is a sequence $\{t_i\}_{i\in \mathbb N}$ with
\[\vint_{B(z, t_{2k})\cap\Omega} |u|\, \rho d\mu\geq \frac23,\ \ \ \vint_{B(z, t_{2k+1})\cap \Omega} |u| \, \rho d\mu \leq \frac13\]
for all $z\in \bx$. We omit the details here but refer to \cite[Remark 3.6]{KNW} and \cite[Lemma 3.6]{KNW-agms} for details of a similar  modification.

Let us go back to find the function $f$ satisfying relation \eqref{eq-example}. Let 
$$h(t)={f(t)}\cdot{t \omega(t)^{1/(\theta-1)}}.$$
Then to find a function $f$ satisfying relation \eqref{eq-example} is equivalent to find a function $h$ satisfying
\begin{equation}\label{eq-example-}\left\{
\begin{array}l
\int_{0}^{1} h(t) \cdot{t^{-1} \omega(t)^{-1/(\theta-1)}}dt=+\infty\\
\int_{0}^{1} h(t)^\theta\cdot t^{-1}\omega(t)^{-1/(\theta-1)}\, dt<+\infty.
\end{array}\right.
\end{equation}
Consider the metric measure space $((0, 1), d_E, \mu_I)$ with $d_E$ being the Euclidean metric and $d\mu_I={t^{-1} \omega(t)^{-1/(\theta-1)}}dt$. Since for any $t_0\in (0, 1)$, we have 
$$\mu_I([t_0, 1))=\int_{t_0}^1 t^{-1} w(t)^{-1/(p-1)}\, dt<\infty.$$ $\mu_I$ is a $\sigma$-finite measure. Moreover, since the relation \eqref{infty} implies that $\mu_I((0, 1))=\infty$, it follows from Lemma \ref{lemma3.1} that $L^\theta((0, 1), \mu_I)\nsubseteq L^1((0, 1), \mu_I)$, i.e., there exists a function $h\colon (0, 1) \rightarrow [0, \infty)$ such that $h\in L^\theta((0, 1), \mu_I)$ but $h\notin L^1((0, 1), \mu_I)$. Such an $h$ satisfies \eqref{eq-example-}.

In conclusion, we have constructed a function $u\in N^{1, \theta}(\Omega, \rho d\mu)$ such that $Tu(z)$ does not exist for any $z\in \partial\Omega$. 
\end{proof}
\end{example}

\begin{rem}
\begin{itemize}
	\item The combination of Theorem \ref{thm-Banach-John} and Example \ref{exam:sharpness} actually shows that the  condition \eqref{eq:John}    characterizes the existence of traces of Sobolev spaces $N^{1, \theta}(\Omega, V, \rho d\mu)$. This characterization is new even when the Banach space $V=\real$, extending and improving the corresponding result of \cite{m17}.
	
	\item Theorem \ref{thm-Banach-John} and  Example \ref{exam:sharpness} are inspired by the recent works from \cite{HWWX} and \cite{KNW}. A full characterization of the existence of traces on regular trees was given in \cite{KNW}. In \cite{HWWX}, the upper half space $\real^{d+1}_+$ with measure $d\mu=|x_{d+1}|^{\theta-1}dx$ and weight function $\varphi(t)=\log^\lambda(4/t)$ was considered. It was shown that the traces of the weighted Sobolev space exist if and only if $\lambda>\theta-1$ for $\theta>1$, which coincides with the relation \eqref{eq:John}. We refer the interested readers to \cite[Theorem 1.2 and Example 1.1]{HWWX} as a special case to understand the more general Theorem \ref{thm-Banach-John} and Example \ref{exam:sharpness}.
\end{itemize}
\end{rem}

\bigskip
\textbf{Acknowledgement.} C.-Y. Guo is supported by the Qilu funding of Shandong University (No. 62550089963197). M. Huang is supported  by NNSF of China (No.11822105). H. Xu is supported by the postdoctor foundation at Shandong University (No. 10000072110302) and the Qilu funding of Shandong University (No. 62550089963197).

%\newpage

\medskip

\noindent Chang-Yu Guo,

\noindent
Research Center for Mathematics and Interdisciplinary
Sciences, Shandong University, Qingdao, Shandong 266237,
P. R. China.

\noindent{\it E-mail address}:  \texttt{changyu.guo@sdu.edu.cn}
\bigskip

\noindent Manzi Huang,

\noindent
MOE-LCSM, School of Mathematics and Statistics, Hunan Normal University, Changsha, Hunan 410081, People's Republic of
China, and School of Mathematical Science, Qufu Normal University, Qufu, Shangdong 273165,
People's Republic of China

\noindent{\it E-mail address}:  \texttt{mzhuang@hunnu.edu.cn}
\bigskip

%\noindent Xiantao Wang,
%
%\noindent
%MOE-LCSM, School of Mathematics and Statistics, Hunan Normal University, Changsha, Hunan 410081, People's Republic of
%China.
%
%\noindent{\it E-mail address}:  \texttt{xtwang@hunnu.edu.cn}
\bigskip

\noindent Zhuang Wang,

\noindent
MOE-LCSM, School of Mathematics and Statistics, Hunan Normal University, Changsha, Hunan 410081, P. R. China.

\noindent{\it E-mail address}:  \texttt{zhuang.z.wang@foxmail.com}, \texttt{zwang@hunnu.edu.cn}
\bigskip

\noindent Haiqing Xu,

\noindent
Research Center for Mathematics and Interdisciplinary
Sciences, Shandong University, Qingdao, Shandong 266237,
P. R. China.

\noindent{\it E-mail address}:  \texttt{hqxu@mail.ustc.edu.cn}
\end{document}